\documentclass[11pt]{article}
\usepackage{graphicx,epsfig}
\usepackage{amsmath,amssymb,euscript}
\usepackage{latexsym}
\usepackage{color}
\newtheorem{theorem}{\bf Theorem}[section]

\def\NN{{\mathbb N}}
\def\RR{{\mathbb R}}

\def\PP{{\mathbb P}}
\def\SS{{\mathbb S}}
\def\C{{\mathcal{C}}}

\def\S{{\mathcal{S}}}

\def\x{{\bf x}}
\def\y{{\bf y}}

\let\OLDthebibliography\thebibliography
\renewcommand\thebibliography[1]{
  \OLDthebibliography{#1}
  \setlength{\parskip}{0pt}
  \setlength{\itemsep}{3mm}
}

\begin{document}
\title{Optimal Lebesgue constants for least squares polynomial approximation on the (hyper)sphere}
\author{Woula Themistoclakis\thanks{C.N.R. National
        Research Council of Italy,
        Istituto per le Applicazioni del Calcolo ``Mauro Picone'',  via P. Castellino, 111, 80131 Napoli, Italy.
        woula.themistoclakis@cnr.it. Partially supported by GNCS-INDAM.}
        \and Marc Van Barel\thanks{KU Leuven, Department of Computer Science, KU Leuven,
Celestijnenlaan 200A,
B-3001 Leuven (Heverlee), Belgium. marc.vanbarel@cs.kuleuven.be.
Supported by
the Research Council KU Leuven,
C1-project (Numerical Linear Algebra and Polynomial Computations),
and by
the Fund for Scientific Research--Flanders (Belgium), EOS Project no 30468160.
}}
\maketitle
\begin{abstract}
We investigate the uniform approximation provided by least squares polynomials on the unit Euclidean sphere $\SS^q$ in $\RR^{q+1}$, with $q\ge 2$.
Like any other polynomial projection, the study concerns the growth, as the degree $n$ tends to infinity, of the associated Lebesgue constant, i.e., of the uniform norm of the least squares operator.  If the least squares polynomial of degree $n$ is based on a set of points, which are nodes of a positive weighted quadrature rule of degree of exactness $2n$, then we state two different sufficient conditions for having an optimal Lebesgue constant that increases with $n$ at the minimal projections order. Hence, under our assumptions least squares and hyperinterpolation polynomials provide a comparable approximation with respect to the uniform norm.
\\[0.2cm]
{\bf keywords:} polynomial approximation on the (hyper)sphere, least squares polynomial, hyperinterpolation, uniform approximation, Lebesgue constant, tensor product quadrature rules.
\\[0.2cm]
\end{abstract}
\section{Introduction}
In recent decades, the polynomial approximation on the sphere
\[
\SS^q:=\left\{\x=(x_0,\ldots,x_{q})\in\RR^{q+1}\ : \sum_{i=0}^q x_i^2=1\right\},\qquad q\ge 2,
 \]
by using function values at a discrete point set $X_N:=\{\xi_1,\ldots,\xi_N\}\subset \SS^q$, has received more and more interest by many authors motivated by the wide field of applications in geophysics, biology and engineering (see, e.g., \cite{b509, W-geoBook, W-GiaMha, W-ReBook, W-GiaSlo} and the references therein).

Limiting our concern to polynomial projections, besides the classical Lagrange interpolation \cite{W-Xu}, we recall the hyperinterpolation polynomials firstly introduced by Sloan in \cite{W-SloFirst} supposing that the point set $X_N$ consists of nodes of a positive weighted quadrature rule of suitable degree of precision. By means of this quadrature rule, hyperinterpolants approximate  Fourier orthogonal projections  w.r.t. the scalar product
\begin{equation}\label{prod}
<f,g> :=\int_{\SS^2}f(\x)g(\x)d\sigma(\x),
\end{equation}
where $d\sigma$ denotes the usual surface measure on $\SS^q$.

It is known \cite{W-Re, W4}  that hyperinterpolation polynomials provide an optimal approximation w.r.t. the uniform norm, but for their construction we need to explicitly know the quadrature weights. This is not necessary if we consider the least squares polynomials, defined as the orthogonal projections w.r.t. the discrete scalar product
\begin{equation}\label{prod-dis}
<f,g>_N :=\sum_{i=1}^Nf(\xi_i)g(\xi_i).
\end{equation}
In the case $q=2$, it has been proved by the authors \cite[Th. 2.3]{W-Marc1} that similarly to hyperinterpolation, also least squares projections have optimal Lebesgue constants w.r.t. the uniform norm, provided that the nodes $\{\xi_j\}$ support a quadrature rule with positive weights (required for hyperinterpolation too) and they are well separated on the sphere.

In this paper we are going to extend this result to the hypersphere case $q\ge 2$ (cf. Theorems \ref{th-LS} and \ref{th-LSequi}).

Moreover, when $q=2$ we focus on the special case of tensor product Gauss--Legendre quadrature rules nodes \cite[Example 6.1]{W4}. These nodes do not satisfy the assumption to be well separated on the sphere. Nevertheless, from our numerical experiments an optimal behavior of the associated Lebesgue constants comes out.

This is justified by a second theorem (cf. Theorem \ref{th-tensor}) that we state in $\SS^q$ with $q\ge 2$, where the assumption on the well separated nodes is replaced by an hypothesis on the quadrature weights, which is certainly satisfied by the tensor product Gaussian quadrature rules.

In the next section we briefly recall some basic properties of Fourier and hyperinterpolation projections. The main theorems are given in Section 3, where some numerical experiments are also given. The proofs are left to Section 4 and Section 5 summarizes the obtained results.
\section{Basic properties of Fourier and hyperinterpolation projections}
Let $\PP_n$ be the space of all spherical polynomials (i.e., polynomials of $q+1$ variables restricted to the sphere $\SS^q$) of degree at most $n$.
It is well--known (see, e.g., \cite{b509}) that
\[
\dim \PP_n= \frac{(2n+q)\Gamma(n+q)}{\Gamma(q+1)\Gamma(n+1)}=:d_n
 \]
and spherical harmonics (i.e., harmonic homogeneous polynomials restricted to $\SS^q$) of degree at most $n$ provide a basis of $\PP_n$, which is orthonormal w.r.t. the scalar product (\ref{prod}).

Moreover, spherical harmonics are related to ultraspherical polynomials of index $q/2-1$ by an addition formula (cf. \cite[(1.6.7)]{b509}), which allows us to write the associated Fourier orthogonal projection $\S_n:L^2(\SS^q)\rightarrow \PP_n$ as follows
\begin{equation}\label{Fourier}
\S_nf(\x)= \frac 1{|\SS^{q-1}|}\int_{\SS^q}K_n(\x\cdot\y)f(\y)d\sigma(\y),\qquad \x\in\SS^q,
\end{equation}
where $|\SS^{q-1}|$ is the surface area of $\SS^{q-1}$, $\x\cdot\y$ denotes the Euclidean scalar product in $\RR^{q+1}$, and
\begin{equation}\label{Darboux}
K_n(t):=K_n(t,1),\qquad t\in [-1,1],
\end{equation}
is the $n$--th Darboux kernel related to the weight function $w(x)=(1-x^2)^{\frac q2-1}$, as defined in \cite{b210}.

Fourier projection, as any other projection onto $\PP_n$, satisfies for all functions $f$ s.t. $\|f\|_\infty:=\sup_{\x\in\SS^q}|f(\x)|<\infty$, the following error estimate
\begin{equation}\label{err}
E_n(f)\le\|f-\S_n f\|_\infty\le \left(1+\|\S_n\|_\infty\right) E_n(f),
\end{equation}
where $E_n(f)$ is the error of best polynomial approximation w.r.t.\ the uniform  norm, i.e.,
\[
E_n(f):=\inf_{P\in\PP_n}\|f-P\|_\infty,
\]
and $\|\S_n\|_\infty$ denotes the so--called Lebesgue constant of $\S_n$, given by
\begin{equation}\label{Leb-Fou}
\|\S_n\|_\infty=\frac 1{|\SS^{q-1}|}\sup_{\x\in\SS^q}\int_{\SS^q}|K_n(\x\cdot\y)|d\sigma(\y).
\end{equation}
 More generally, we recall that the Lebesgue constant of any projection $T_n$ is defined as the following operator norm
\[
\|T_n\|_\infty:=\sup_{\|f\|_\infty\le 1}\|T_nf\|_\infty,
\]
and its behaviour as $n\rightarrow + \infty$ strongly influences the quality of the approximation.

It is known (see, e.g., \cite{b509, W-Da, W-ReBook}) that the previous Fourier projection $\S_n$ is the projection onto $\PP_n$ having minimal Lebesgue constant. More precisely, if we denote by ${\cal T}_n$ the class of all the polynomial projections onto $\PP_n$, then for sufficiently large $n$, we have
\begin{equation}\label{min-norm}
\|T_n\|_\infty\ge \|\S_n\|_\infty\sim n^\frac{q-1}2,\qquad \forall T_n\in{\cal T}_n,
\end{equation}
where throughout the paper by $a_n\sim b_n$ we mean that $c_1 a_n\le b_n\le c_2 a_n$ being $c_1,c_2>0$ independent of $n$.

However, the approximation $\S_nf$ requires the computation of the Fourier coefficients that are integrals of the unknown function $f$.
If we suppose to know $f$ only at a discrete point set $X_N:=\{\xi_1,\ldots,\xi_N\}$ such that the quadrature rule
\begin{equation}\label{quad-2n}
\int_{\SS^q} f(\x)d\sigma(\x)=\sum_{i=1}^N\lambda_i f(\xi_i),\qquad \lambda_i>0,\qquad\quad\forall f\in\PP_{2n},
\end{equation}
holds, then we can discretize $\S_nf$ by applying (\ref{quad-2n}) to (\ref{Fourier}). In this way, we get the following polynomial of degree at most $n$ \cite{W-SloFirst, W4}
\begin{equation}\label{hyper}
L_{n}f(\x)=\frac 1{|\SS^{q-1}|}\sum_{i=1}^N\lambda_{i}f(\xi_i)K_{n}(\xi_i\cdot \x),\qquad \x\in\SS^q,
\end{equation}
which is usually called {\it hyperinterpolation polynomial} because it is based on the function values at a number of nodes $N$ that is greater than $d_n$, the dimension of $\PP_n$ \cite{W-Ba, W-Re}.

The double degree of exactness in (\ref{quad-2n}) assures that  $L_{n}$ is a discrete polynomial projection onto $\PP_{n}$, namely
\begin{equation}\label{inva-hyper}
L_nP=P,\qquad \forall P\in\PP_n.
  \end{equation}
  Moreover, it is known that  the Lebesgue constants $\|L_{n}\|_\infty$ increase with $n$ at the order of the minimal projections, i.e., for all sufficiently large $n\in\NN$, we have \cite{W-Re, W4}
\begin{equation}\label{hyper-norm}
\|L_{n}\|_\infty\sim\|\S_n\|_\infty\sim n^\frac{q-1}2.
\end{equation}
\section{On least squares polynomial approximation}
A different kind of discrete polynomial projection is given by the least squares approximations $\tilde\S_nf\in \PP_n$, defined by
\begin{equation}\label{LS-min}
\sum_{i=1}^N[f(\xi_i)-\tilde\S_nf(\xi_i)]^2=\min_{P\in\PP_n}
\sum_{i=1}^N[f(\xi_i)-P(\xi_i)]^2.
\end{equation}
In explicit form, for all $\x\in\SS^q$, the least squares polynomial $\tilde\S_nf(\x)$ related to the point set $X_N=\{\xi_1,\ldots,\xi_N\}\subset\SS^q$ is given by
\begin{equation}\label{LS-sum}
\tilde\S_{n}f(\x)=\sum_{i=1}^Nf(\xi_i)H_{n}(\x, \xi_i),\qquad H_n(\x,\y):=\sum_{r=1}^{d_n}I_r(\x)I_r(\y),
\end{equation}
where $\{I_r : r=1,\ldots,d_n\}$ is a basis of $\PP_n$ orthonormal w.r.t. the discrete scalar product
defined in (\ref{prod-dis}).
Moreover, we observe that
\begin{equation}\label{inva}
P(\x)=\sum_{i=1}^NP(\xi_i)H_{n}(\x, \xi_i),\qquad  \forall P\in\PP_n, \qquad \forall\x\in\SS^q .
\end{equation}
With respect to the hyperinterpolation $L_nf$, the least squares polynomial $\tilde\S_nf$ does not require to know any quadrature weight, neither any quadrature rule is indeed necessary for its definition.

Concerning the Lebesgue constant $\|\tilde\S_n\|_\infty$, for the $2$--sphere case (i.e., $q=2$) in \cite{W-Marc1} it has been proved that $\|\tilde\S_n\|_\infty\sim \|L_n\|_\infty$ holds if the point set $X_N=\{\xi_1,\xi_2,\ldots,\xi_N\}\subset \SS^2$ is such to support a positive weighted quadrature rule of degree of exactness $2n$ and if  the following Marcinkievicz type inequality holds
\[
\frac 1{n^2}\sum_{i=1}^N |Q(\xi_i)|\le\C \|Q\|_{L^1(\SS^2)},\qquad \forall Q\in\PP_{n},\qquad \C\ne\C(n,N,Q),
\]
where throughout the paper we denote by $\C$ a positive constant, which can take different values at the different occurrences, and we write $\C\ne\C(n,N,Q,..)$ to mean that $\C$ is independent of $n,N,Q,...$

The next theorem generalizes \cite[Th. 2.3]{W-Marc1} to any dimension $q\ge 2$.
\begin{theorem}\label{th-LS}
Let the point set $X_N=\{\xi_1,\xi_2,\ldots,\xi_N\}\subset \SS^q$   and $n\in\NN$ be such that (\ref{quad-2n}) holds. 
Moreover, suppose that
\begin{equation}\label{Marci-q}
\frac 1{n^q}\sum_{i=1}^N |Q(\xi_i)|\le\C \|Q\|_{L^1(\SS^q)},\qquad \forall Q\in\PP_{n},\qquad \C\ne\C(n,N,Q).
\end{equation}
Then for all sufficiently large $n\in\NN$, the Lebesgue constant of the least squares polynomial of degree $n$ associated to $X_N$, satisfies
\begin{equation}\label{Leb-LS}
\|\tilde\S_n\|_\infty\sim n^{\frac{q-1}2}.
\end{equation}
\end{theorem}
Let $card(A)$ denote the cardinality of the set $A$ and let $d(\x,\y):=\arccos(\x\cdot\y)$ be the geodesic distance of $\x,\y\in\SS^q$. In \cite[Th. 2.1]{W-Dai} it has been proved that
\begin{equation}\label{hp-sep}
\sup_{\x\in\SS^q} card\left(\left\{\xi_i\in X_N : d(\xi_i, \x)\le \frac 1n\right\}\right)\le \C, \qquad \C\ne\C(n,N),
\end{equation}
is a necessary and sufficient condition for (\ref{Marci-q}), so that the previous theorem is equivalent to the following
\begin{theorem}\label{th-LSequi}
Let the point set $X_N=\{\xi_1,\xi_2,\ldots,\xi_N\}\subset \SS^q$ be such that (\ref{hp-sep}) and (\ref{quad-2n}) holds.  Then for sufficiently large $n\in\NN$, we have
\[
\|\tilde\S_n\|_\infty\sim n^{\frac{q-1}2}.
\]
\end{theorem}
The assumption (\ref{hp-sep}) is also required to state the existence of positive weighted quadrature rules (see, e.g., \cite{W-Dai, r938, W-Po}). Nevertheless there exist several positive quadrature rules not satisfying (\ref{hp-sep}).
This is the case of tensor product Gauss--Legendre quadrature rules deduced in \cite[Example 6.1]{W4} for $q=2$, by combining the trigonometric rectangular rule (exact for trigonometric polynomials of degree $\le 2n+1$)
\[
\int_0^{2\pi}g(\Phi)d\Phi=\frac\pi{n+1}\sum_{k=0}^{2n+1}g\left(\Phi_k\right),\qquad \Phi_k:=\frac{k\pi}{n+1},
\]
and the $(n+1)$--point Gauss--Legendre quadrature rule
\begin{equation}\label{GL-quadrule}
\int_{-1}^1G(z)dz=\sum_{j=1}^{n+1}\nu_jG(z_j). 
\end{equation}
The resulting tensor product rule has degree of precision $2n+1$ and it is based on $N=2(n+1)^2$ points. It looks like
\begin{equation}\label{tensor}
\int_{\SS^2}f(\x)d\sigma(\x)=\sum_{k=0}^{2n+1}\sum_{j=1}^{n+1} \frac {\pi \nu_j}{n+1} f(\xi_{j,k}),
\qquad \forall f\in\PP_{2n+1},
\end{equation}
where each node $\xi_{j,k}$ has azimuthal angle $\Phi_k$ and polar angle $\theta_j=\arccos z_j$.

The main advantage of tensor product rules is the explicit knowledge of the quadrature weights and nodes, but the latter have the disadvantage of not being well--separated on the sphere.

\begin{figure}
\includegraphics[scale=0.5]{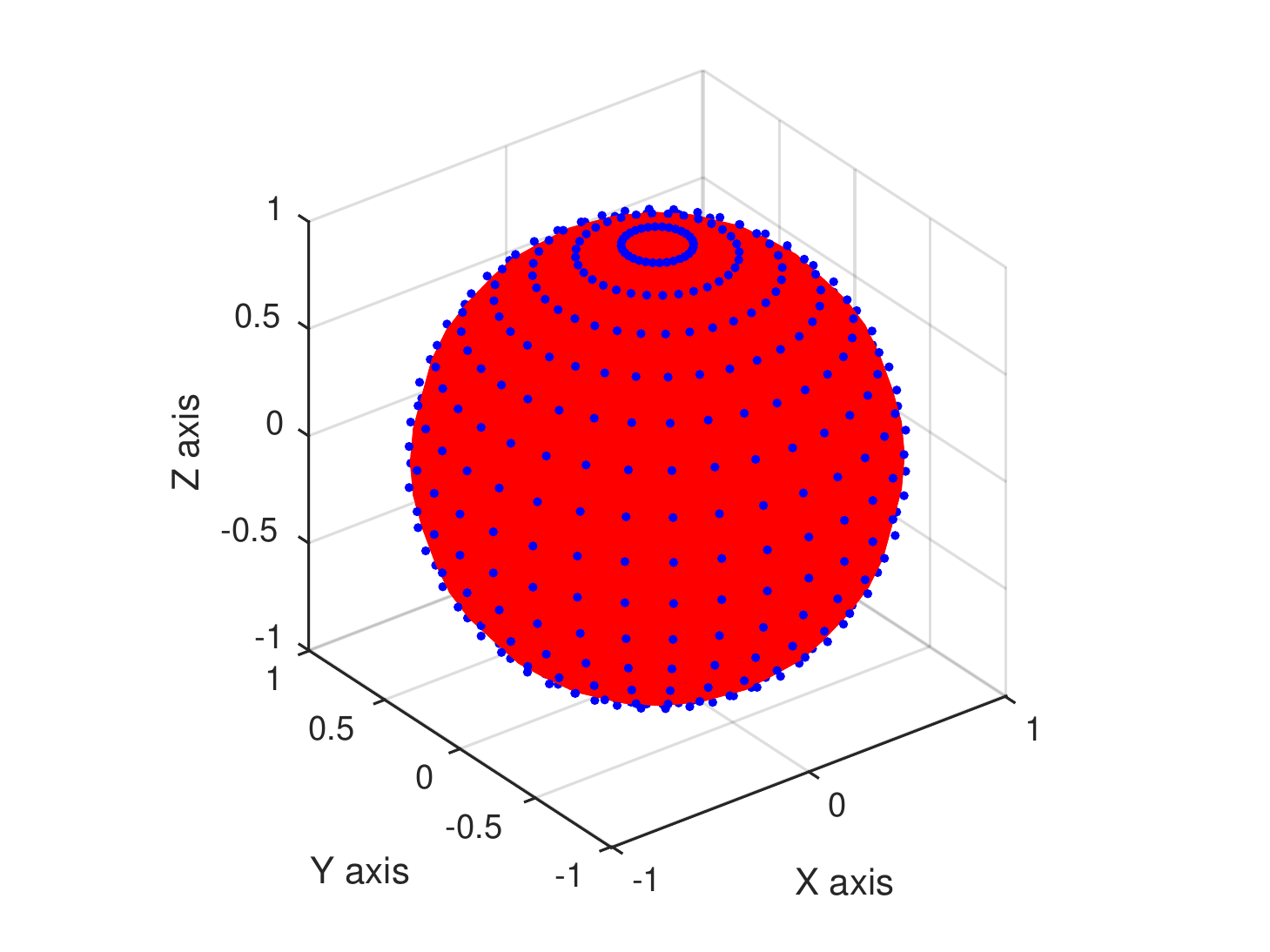}
\includegraphics[scale=0.5]{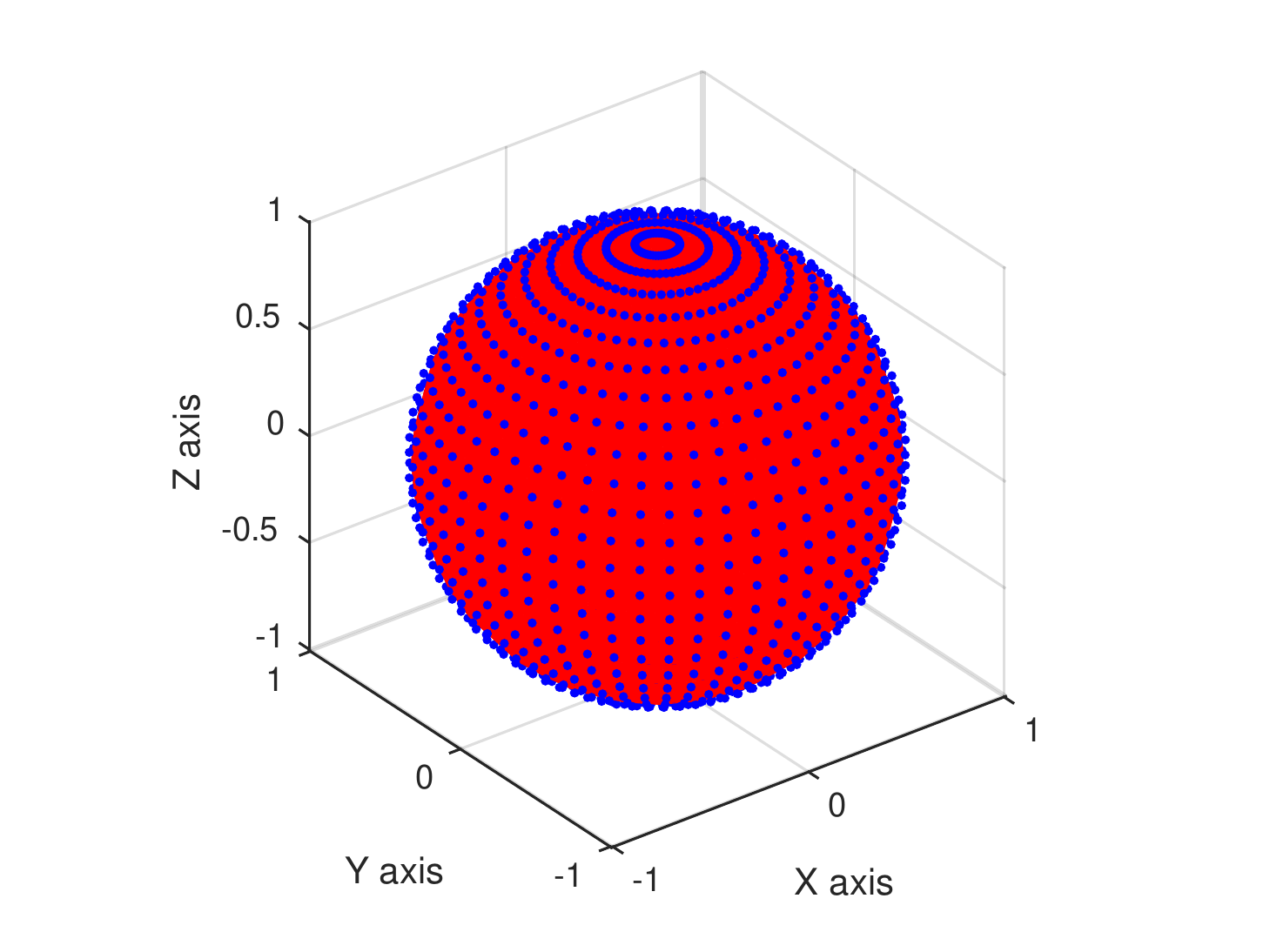}
\caption{Examples of the tensor product Gauss--Legendre quadrature nodes related to degrees of precision $31$ and $51$, i.e., having $N=512$ (left) and $N=1352$ (right) points.\label{fig001}}
\end{figure}

This can be seen by Figure~\ref{fig001}, which shows how the nodes
\begin{equation}\label{tensor-pt}
\tilde X_N:=\{\xi_{j,k} : k=0,...,2n+1, \ j=1,...,n+1\}
\end{equation}
are distributed on the sphere $\SS^2$ for degree of precision $31$, i.e.,  $n=15$, $N=512$, and degree
of precision $51$, i.e., $n=25$, $N = 1352$.

We recall (see, e.g., \cite{HardMichSaff2016}) that a measure of the uniformity of a sampling set $X_N=\{\xi_1,\ldots,\xi_N\}$ is given by the mesh norm $\delta_{X_N}$ and separation distance $\gamma_{X_N}$ defined by
\begin{eqnarray}\label{delta}
\delta_{X_N}&:=&\max_{\x\in\SS^2}\min_{1\le i\le N} d(\x,\xi_i),\\
\label{gamma}
\gamma_{X_N}&:=&\min_{ i\ne j}d(\xi_i,\xi_j),
\end{eqnarray}
and a sequence of point configurations $\{X_N\}_N$ is said to be quasi--uniform if the mesh ratio $\delta_{X_N}/\gamma_{X_N}$ is bounded as $N\rightarrow +\infty$.

\begin{figure}
\includegraphics[scale=0.5]{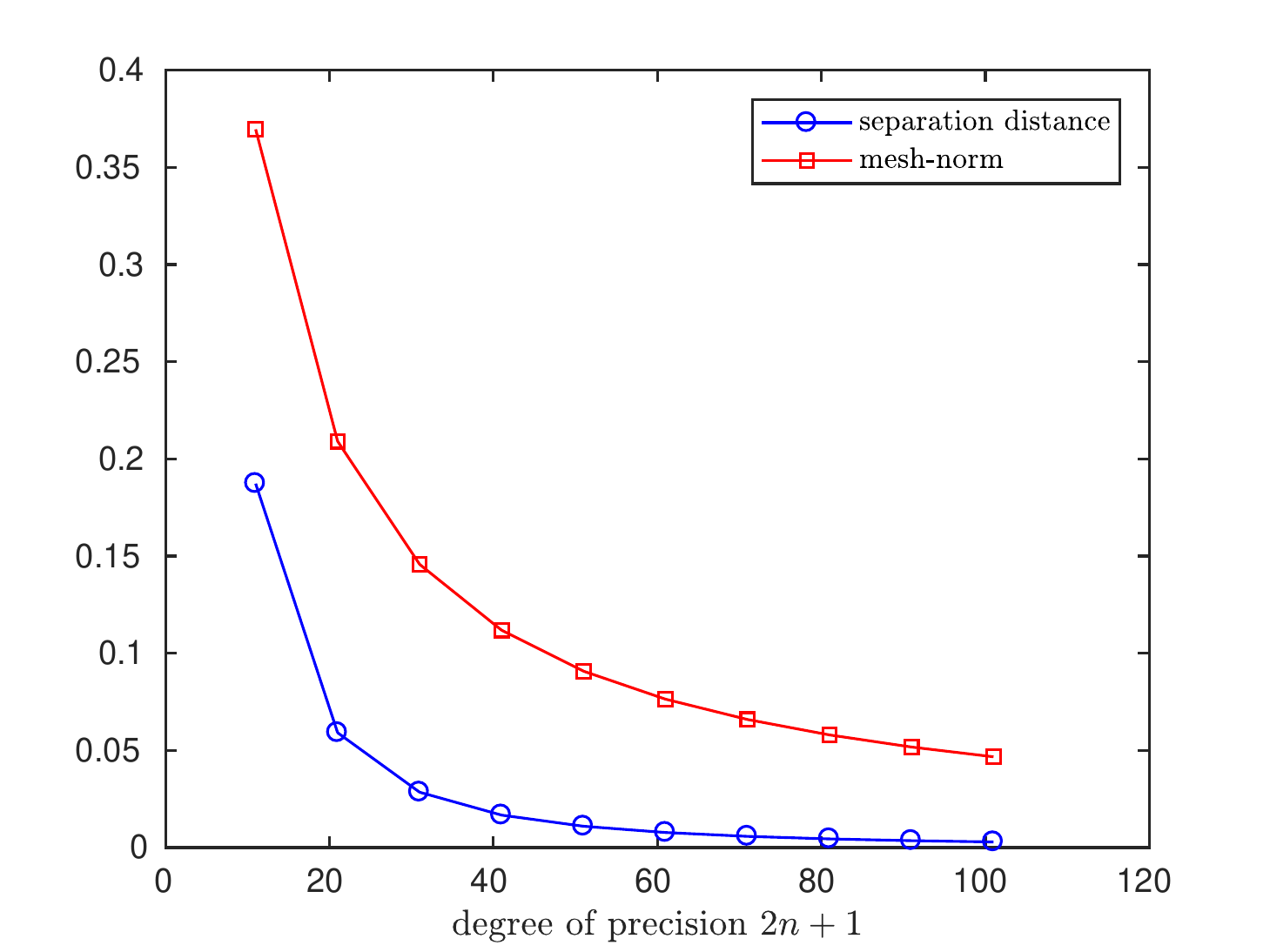}
\includegraphics[scale=0.5]{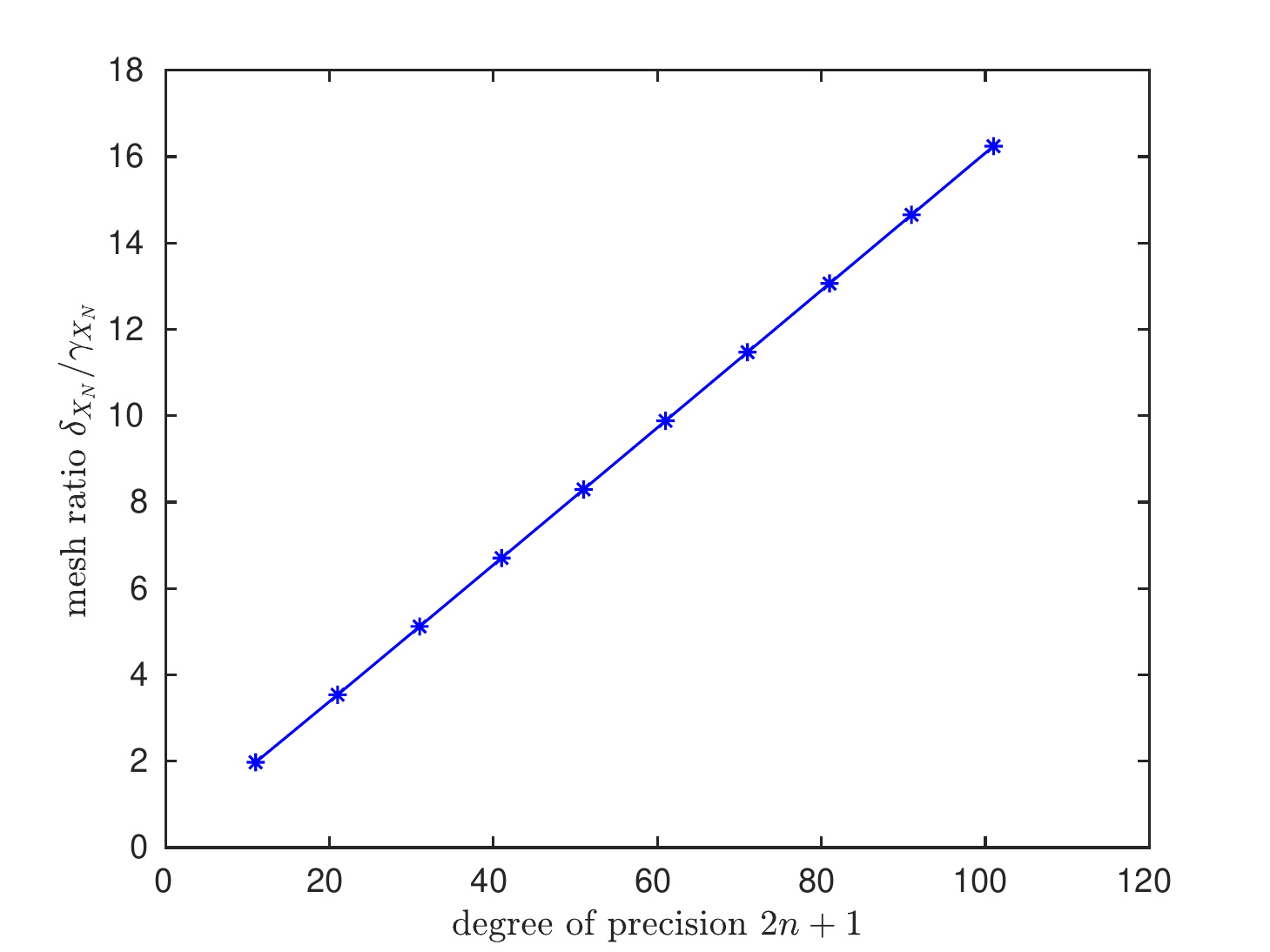}
\caption{The values of the separation distance $\gamma_{\tilde X_N}$ and the mesh-norm $\delta_{\tilde X_N}$ (left) and the mesh ratio $\delta_{\tilde X_N} / \gamma_{\tilde X_N}$ (right) for point sets $\tilde X_N$ with $N = 2(n+1)^2$ for $n = 5, 10, 15,\ldots,50$.
\label{fig011}}
\end{figure}

Figure~\ref{fig011} displays some values of the mesh norm, separation distance and mesh ratio for the point set $\tilde X_N$ in (\ref{tensor-pt}).
To estimate the mesh norm $\delta_{X_N}$, instead of taking the maximum over  the set of all
points of the sphere, the maximum is computed over a point set with a number of points considerably larger than the number of points for which we want
to approximate the mesh norm. To this end we consider the ``spiral points'' as defined in  \cite{RakhSaffZhou1995} (see also the overview paper \cite{HardMichSaff2016}).  In the sequel we refer to this point set as the point set of second type.
These can be computed very efficiently and seem to be uniformly distributed
over the unit sphere where each of the points seems to be well separated from the others.
To estimate $\delta_{\tilde X_N}$ we considered a point set of the second type having $16$ times more points.

It is evident that the distribution of the nodes in $\tilde X_N$ is not uniform and indeed it turns out that
\[
\delta_{\tilde X_N}\le\frac\C n, \qquad \mbox{but}\qquad \gamma_{\tilde X_N}\ge\frac \C{n^2},\qquad\qquad \C\ne\C(n,N).
\]
Hence, we can say that the previous tensor product nodes provide an optimal hyperinterpolation polynomial, but from a theoretical point of view, up to now nothing can be said regarding the least squares polynomial, since the assumption (\ref{hp-sep}) of the previous theorem is not satisfied.

\begin{figure}
\includegraphics[scale=0.85]{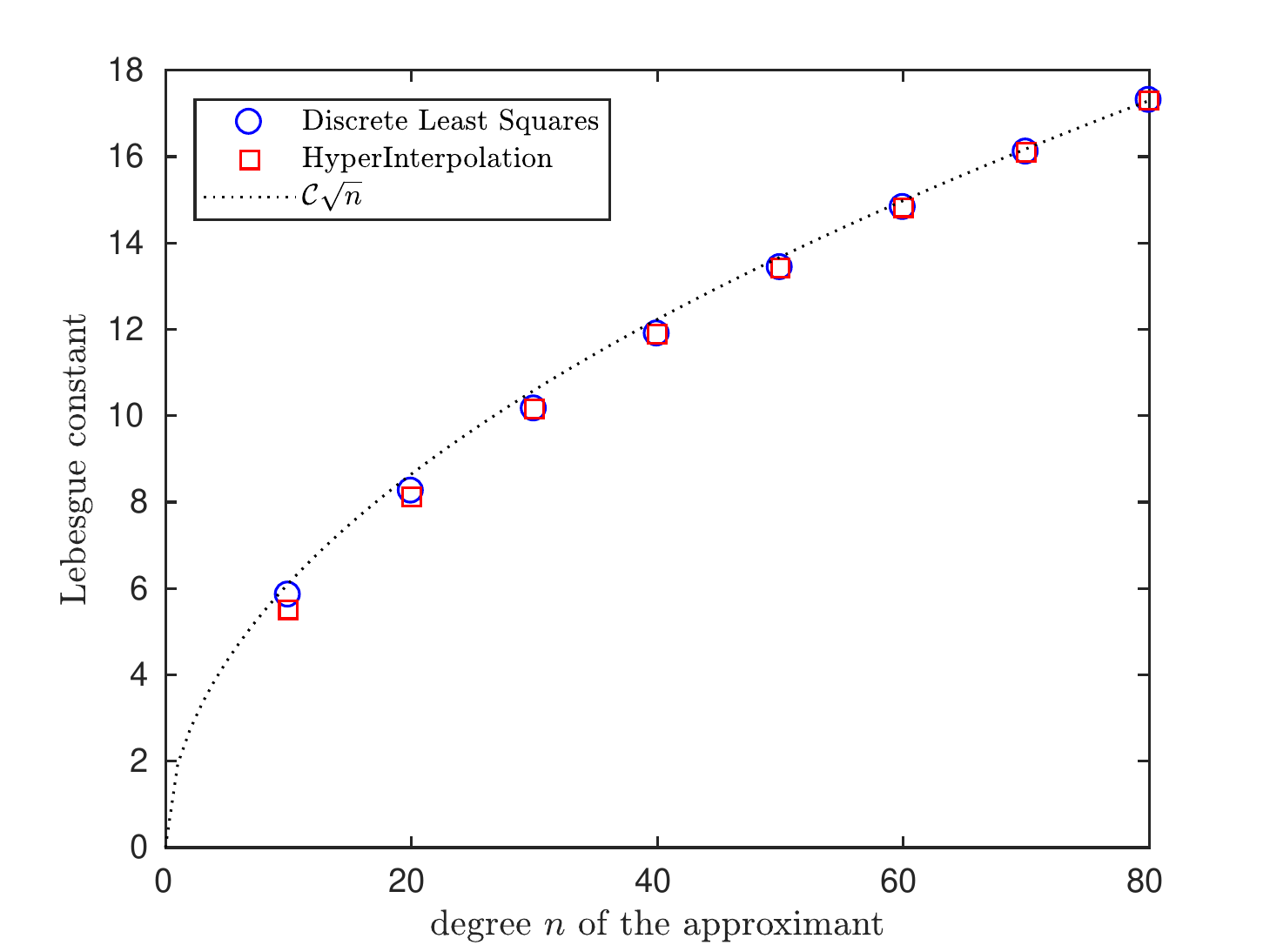}
\caption{The values of the Lebesgue constants of the least squares operator $\tilde\S_n$ and the hyperinterpolation operator $L_n$ for the degrees $n=10,20,\ldots,80$ with corresponding point set $\tilde X_N$ having $N = 2(n+1)^2$ points.\label{fig021}}
\end{figure}

Now we investigate numerically the behaviour of the Lebesgue constants of both least squares and hyperinterpolation polynomials of degree $n$ related to the previous point set $\tilde X_N$.
 To this end, we'll estimate the uniform norm of the corresponding operators by taking a larger point set  of the second type containing $4N$   points.
Figure~\ref{fig021} shows the results. The circles and squares indicate the Lebesgue constant for the least squares operator and the hyperinterpolation operator, respectively, when we take for degree $n$ on the horizontal axis
the corresponding point set $\tilde X_N$  related to the degree of exactness $2n+1$, i.e., having $N = 2(n+1)^2$ points.
 The figure shows for $\|\tilde\S_n\|_\infty$ the same optimal behavior as for $\|L_n\|_\infty$, i.e., the $\sqrt{n}$ behaviour as indicated by the dotted line.

In order to explain such a numerical output in the case of the Gauss--Legendre tensor product rule (\ref{tensor}), we recall that the Legendre zeros
 \[
 z_0:=-1<z_1<\ldots <z_{n+1}<1=:z_{n+2}
 \]
 are arcsin distributed on $[-1,1]$ and that for the weights $\nu_i$ of the Gauss-Legendre quadrature rule (\ref{GL-quadrule}), $\nu_i\sim (z_{i+1}-z_i)$ holds uniformly w.r.t.\ $i$ and $n$ \cite{W-Ne}. Consequently
\[
\nu_i\sim \nu_{i+ 1},\qquad i=1,\ldots, n,
\]
holds uniformly w.r.t.\ $i$ and $n$. This is indeed the property replacing (\ref{Marci-q}) or (\ref{hp-sep}), in order to get the same result as in the previous theorems.
\begin{theorem}\label{th-tensor}
Let $X_N=\{\xi_1,\ldots,\xi_N\}\subset \SS^q$ and $\lambda_1\ge\lambda_2\ge\ldots\ge\lambda_N>0=:\lambda_{N+1}$ be such that
\[
\int_{\SS^q} f(\x)d\sigma(\x)=\sum_{i=1}^N\lambda_i f(\xi_i),\qquad\forall f\in\PP_{2n},
\]
holds for $n\in\NN$.
Moreover, suppose that the quadrature weights are such that
\begin{equation}\label{hp-tensor}
\lambda_i\le \C \lambda_{i+ 1},\qquad i=1,\ldots, N-1, \qquad \C\ne \C(n,N,i).
\end{equation}
Then for all sufficiently large $n$, we have
$\displaystyle
\|\tilde\S_n\|_\infty\sim n^\frac{q-1}2.
$
\end{theorem}
In contrast to the previous two theorems, Theorem \ref{th-tensor} results to be applicable to all tensor product rules that usually have the nodes very close to each other near the poles, so (\ref{hp-sep}) as well as (\ref{Marci-q}) do not generally hold, but (\ref{hp-tensor}) holds.
\section{Proofs}
\subsection{Proof of Theorem \ref{th-LS}.}
From (\ref{LS-sum}) we deduce that
\[
\|\tilde\S_n\|_\infty=\sup_{\x\in\SS^q} \left[\sum_{k=1}^N|H_{n}(\x, \xi_k)|\right].
\]
Hence, due to (\ref{min-norm}), it is sufficient to prove that
\begin{equation}\label{eq-th}
\sum_{k=1}^N|H_n(\x,\xi_k)|\le\C n^{\frac{q-1}2}, \qquad \forall\x\in\SS^q, \qquad \C\ne\C(n,N,\x).
\end{equation}
Let us first prove (\ref{eq-th}) when $\x\in X_N$.

To this aim we observe that $P(\xi_k):=H_n(\x,\xi_k)$ is a spherical polynomial of degree $n$ w.r.t.\ the variable $\xi_k$. Consequently, recalling that $L_nP=P$, we get
\begin{equation}\label{trans1}
H_n(\x,\xi_k)=\frac 1{|\SS^{q-1}|}\sum_{i=1}^N\lambda_i  H_n(\x,\xi_i)
K_n(\xi_k\cdot\xi_i),\qquad k=1,\ldots,N,
\end{equation}
where without losing the generality, we assume that $\lambda_i$ are labeled in non increasing order, namely
\[
\lambda_1\ge\lambda_2\ge\ldots \ge\lambda_N>\lambda_{N+1}:=0.
\]
Then, by applying  the following summation by parts formula
\begin{equation}\label{sum-part}
\sum_{i=1}^Na_ib_i=a_N\sum_{i=1}^N b_i+\sum_{i=1}^{N-1}(a_i-a_{i+1})\sum_{j=1}^ib_j,
\end{equation}
we get
\begin{equation}\label{eq-start1}
H_n(\x,\xi_k)= \frac 1{|\SS^{q-1}|}\sum_{i=1}^{N}(\lambda_i-\lambda_{i+1})\sum_{j=1}^iH_n(\x,\xi_j)
K_n(\xi_k\cdot\xi_j).
\end{equation}
Consequently, by taking into account that $\lambda_i-\lambda_{i+1}\ge 0$, we have
\begin{eqnarray*}
\sum_{k=1}^N|H_n(\x,\xi_k)|&\le&\frac 1{|\SS^{q-1}|} \sum_{k=1}^N\sum_{i=1}^{N}(\lambda_i-\lambda_{i+1})
\left|\sum_{j=1}^iH_n(\x,\xi_j)K_n(\xi_k\cdot\xi_j)\right|\\
&=&\frac 1{|\SS^{q-1}|}
\sum_{i=1}^{N}(\lambda_i-\lambda_{i+1})\sum_{k=1}^N\left|\sum_{j=1}^iH_n(\x,\xi_j)
K_n(\xi_k\cdot\xi_j)\right|\\
&\le&\frac 1{|\SS^{q-1}|}
\sup_{1\le i\le N}\left(\sum_{k=1}^N\left|\sum_{j=1}^iH_n(\x,\xi_j)K_n(\xi_k\cdot\xi_j)\right|\right)
\sum_{i=1}^{N}(\lambda_i-\lambda_{i+1})
\\
&=& \frac{\lambda_1}{|\SS^{q-1}|}\sup_{1\le i\le N}\sum_{k=1}^N
\left|\sum_{j=1}^iH_n(\x,\xi_j)K_n(\xi_k\cdot\xi_j)\right|,
\end{eqnarray*}
and recalling that \cite[Lemma 5.4.3]{b509}
\begin{equation}\label{lambda}
\lambda_i\le \frac\C{n^q},\qquad i=1,\ldots, N,\qquad \C\ne \C(n,N,i),
\end{equation}
we obtain
\[
\sum_{k=1}^N|H_n(\x,\xi_k)|\le \frac \C{n^q}\sup_{1\le i\le N} \sum_{k=1}^N
\left|\sum_{j=1}^iH_n(\x,\xi_j)K_n(\xi_k\cdot\xi_j)\right|,\qquad \C\ne\C(n,N,\x).
\]
Hence, set for any $n\in\NN$ and $\x\in X_N$
\[
A_i:=\frac 1{n^q} \sum_{k=1}^N
\left|\sum_{j=1}^iH_n(\x,\xi_j)K_n(\xi_k\cdot\xi_j)\right|, \qquad i=1,\ldots,N,
\]
to get the statement when $\x\in X_N$, we are going to prove that
\begin{equation}\label{tesi}
A_i\le\C n^{\frac{q-1}2}, \qquad i=1,\ldots,N,\qquad \C\ne\C(n,N,\x),
\end{equation}
holds for all sufficiently large $n\in\NN$ and any $\x\in X_N$.

For the case $i=N$, note that by (\ref{inva}) we get
\[
A_N:=\frac 1{n^q} \sum_{k=1}^N
\left|\sum_{j=1}^NH_n(\x,\xi_j)K_n(\xi_k\cdot\xi_j)\right|= \frac 1{n^q} \sum_{k=1}^N
\left|K_n(\xi_k\cdot\x)\right|.
\]
Moreover, (\ref{Marci-q}) and (\ref{Leb-Fou})  imply
\[
\frac 1{n^q} \sum_{k=1}^N
\left|K_n(\xi_k\cdot\xi)\right|\le \C\int_{\SS^q}\left|K_n(\y\cdot\xi)\right|d\sigma(\y)\le\C \|\S_n\|_\infty,\qquad \forall\xi\in\SS^q,
\]
and hence by (\ref{min-norm}) we have
\begin{equation}\label{eq-darboux}
\frac 1{n^q} \sum_{k=1}^N
\left|K_n(\xi_k\cdot\xi)\right|\le \C \|\S_n\|_\infty\le\C n^{\frac{q-1}2},\qquad \forall\xi\in\SS^q,
\qquad \C\ne\C(n,N,\xi).
\end{equation}
So, we conclude that
\begin{equation}\label{AN}
A_N\le\C n^{\frac{q-1}2}, \qquad \forall \x\in X_N,\qquad \C\ne\C(n,N,\x) .
\end{equation}
As regards the case $1\le i<N$,
we observe that for any pair of nodes $\xi_l,\xi_j\in X_N$, we have
\begin{equation}\label{CSW}
|H_n(\xi_l,\xi_j)|\le\sum_{r=1}^{d_n}|I_r(\xi_l)I_r(\xi_j)|\le \left(\sum_{r=1}^{d_n}|I_r(\xi_l)|^2\right)^\frac 12
\left(\sum_{r=1}^{d_n}|I_r(\xi_j)|^2\right)^\frac 12.
\end{equation}
On the other hand, we point out that the existence of (\ref{quad-2n}) implies that $d_n< N$ for sufficiently large $n$ (see, e.g., \cite[p. 274]{W-Re}). Consequently, the matrix consisting of the orthonormal columns $[I_k(\xi_1),\ldots I_k(\xi_N)]^T$, $k=1,\ldots, d_n$, namely the matrix
\[
I := [I_k(\xi_h)]_{h=1,...,N}^{k=1,\ldots,d_n}
 \]
is rectangular, but it can be extended by additional columns to form a square orthogonal matrix
\[
Q=[Q_{h,k}]_{h=1,...,N}^{k=1,...,N}, \qquad \mbox{such that}\qquad Q_{h,k}=I_k(\xi_h),\quad \forall k\le d_n.
 \]
Thus we have
\begin{equation}\label{Q}
\sum_{k=1}^{d_n}|I_k(\xi_j)|^2\le \sum_{k=1}^{N}|Q_{j,k}|^2 =1,\qquad j=1,\ldots,N,
\end{equation}
and assembling (\ref{CSW}) and (\ref{Q}), we conclude that
\begin{equation}\label{ls-ker1}
|H_n(\xi_l,\xi_j)|\le 1,\qquad \forall \xi_l,\xi_j\in X_N.
\end{equation}
By means of (\ref{ls-ker1}), we deduce
\[
A_1:=\frac 1{n^q} \sum_{k=1}^N
\left|H_n(\x,\xi_1)K_n(\xi_k\cdot\xi_1)\right|\le \frac 1{n^q} \sum_{k=1}^N
\left|K_n(\xi_k\cdot\xi_1)\right|,\qquad \forall \x\in X_N,
\]
and using (\ref{eq-darboux}), we get
\begin{equation}\label{A1}
A_1\le\C n^{\frac{q-1}2}, \qquad \forall \x\in X_N,\qquad \C\ne \C(n,N,\x) .
\end{equation}
Similarly, for any $1\le i<N$ and $\x\in X_N$, by (\ref{ls-ker1}) and (\ref{eq-darboux}), we have
\begin{eqnarray}
\nonumber
|A_{i+1}- A_i| &\le&\frac 1{n^q} \sum_{k=1}^N
\left|H_n(\x,\xi_{i+1})K_n(\xi_k\cdot\xi_{i+1})\right|\\
\nonumber
&\le& \frac 1{n^q} \sum_{k=1}^N\left|K_n(\xi_k\cdot\xi_{i+1})\right|\\
\label{Ai}
&\le&\C n^{\frac{q-1}2},\qquad \C\ne\C(n,N,\x, i).
\end{eqnarray}
In conclusion, let us show that (\ref{AN}), (\ref{A1}) and (\ref{Ai}) imply that as $n\rightarrow + \infty$ (\ref{tesi}) holds for all $x\in X_N$.

Indeed, if ad absurdum there exists $\x\in X_N$ s.t. for some index $l$ we have that
\[
\limsup_{n\rightarrow +\infty}\  n^\frac{1-q}2 A_l=+\infty,
\]
then (\ref{AN})  and (\ref{A1}) imply $1<l<N$, and from (\ref{Ai}) we deduce that we also have
\[
\limsup_{n\rightarrow +\infty}\ n^\frac{1-q}2 A_{l\pm 1} =+\infty.
\]
Thus, by iterating the reasoning, we arrive to contradict (\ref{AN}) or (\ref{A1}).

Hence, we conclude that (\ref{eq-th}) holds for all $\x\in X_N$.

For arbitrary $\x\in\SS^q$, we reason analogously, but we start applying the invariance property $L_nP=P$ to the polynomials $P(\x)=H_n(\x,\xi_k)$, with $k=1,\ldots,N$. Hence, instead of (\ref{trans1}) we get
\begin{equation}\label{trans2}
H_n(\x,\xi_k)=\frac 1{|\SS^{q-1}|}\sum_{i=1}^N\lambda_i  H_n(\xi_i,\xi_k)
K_n(\x\cdot\xi_i),\qquad k=1,\ldots,N,
\end{equation}
which differs from (\ref{trans1}) by the exchanged position of the variables $\x$ and $\xi_k$ at the right--hand sides.

Consequently, by using (\ref{sum-part}) and (\ref{lambda}) as before, we deduce
\begin{eqnarray*}
\sum_{k=1}^N|H_n(\x,\xi_k)|&=&\frac 1{|\SS^{q-1}|} \sum_{k=1}^N\left|\sum_{i=1}^N\lambda_i  H_n(\xi_i,\xi_k)
K_n(\x\cdot\xi_i)\right|\\
&\le& \frac\C{n^q}\max_{1\le i\le N}\left(\sum_{k=1}^N\left|\sum_{j=1}^i H_n(\xi_j,\xi_k)
K_n(\x\cdot\xi_j)\right|\right), \qquad \C\ne\C(n,N,\x).
\end{eqnarray*}
Then, for arbitrarily fixed $n\in\NN$ and $\x\in\SS^q$, we set
\[
B_i:=\frac 1{n^q}\sum_{k=1}^N\left|\sum_{j=1}^i H_n(\xi_j,\xi_k)K_n(\x\cdot\xi_j)\right|,\qquad i=1,\ldots,N.
\]
When $i=N$, by virtue of (\ref{inva}) and (\ref{eq-darboux}), we have
\begin{equation}\label{BN}
B_N=\frac 1{n^q}\sum_{k=1}^N\left|K_n(\x\cdot\xi_k)\right|\le\C n^{\frac{q-1}2},\qquad \forall \x\in\SS^q,\qquad \C\ne\C(n,N,\x).
\end{equation}
Moreover, recalling that (see, e.g., \cite{r938, b210})
\begin{equation}\label{sup-Darboux}
|K_n(\x\cdot\y)|\le \sup_{|t|\le 1}|K_n(t)|= K_n(1)\sim n^q,\qquad \forall\x,\y\in\SS^q,
\end{equation}
and taking into account that we have already proved (\ref{eq-th}) in $X_N$, for all $\x\in\SS^q$ we get
\begin{eqnarray}\nonumber
B_1&:=&\frac 1{n^q}\sum_{k=1}^N\left| H_n(\xi_1,\xi_k)K_n(\x\cdot\xi_1)\right|\le
\C\sum_{k=1}^N\left|H_n(\xi_1\cdot\xi_k)\right|\\
\label{B1}
&\le&\C n^{\frac{q-1}2},\qquad\qquad\C\ne\C(n,N,\x,\xi_1),
\end{eqnarray}
as well as, for all $i=1,\ldots,N-1$, and any $\x\in\SS^q$, we deduce
\begin{eqnarray}\nonumber
|B_{i+1}- B_i|&\le&\frac 1{n^q}\sum_{k=1}^N\left|H_n(\xi_{i+1},\xi_k)K_n(\x\cdot\xi_{i+1})\right|\\
\nonumber
&\le& \C\sum_{k=1}^N\left|H_n(\xi_{i+1},\xi_k)\right|\\
\label{Bi}
&\le& \C n^{\frac{q-1}2},\qquad\qquad\C\ne\C(n,N,\x,i).
\end{eqnarray}
In conclusion, similarly to the case $\x\in X_N$,  from (\ref{BN}), (\ref{B1}) and (\ref{Bi}) we deduce that for all sufficiently large $n\in\NN$, and any $\x\in \SS^q$, we have
\[
B_i\le \C n^{\frac{q-1}2},\qquad i=1,\ldots,N,\quad\qquad\C\ne\C(n,N,\x),
\]
and the statement follows in the case $\x\in\SS^q$ too.
\subsection{Proof of Theorem \ref{th-tensor}.}
Following the same reasoning of the previous proof, we arrive to say that it is sufficient to state that for all sufficiently large $n\in\NN$, we have
\[
\sum_{k=1}^N|H_n(\x,\xi_k)|\le\C n^{\frac{q-1}2}, \qquad \forall\x\in X_N ,\qquad \C\ne\C(n,N,\x).
\]
Note that, by using (\ref{eq-start1}) and $\lambda_i-\lambda_{i+1}\ge 0$, we get
\begin{eqnarray*}
|\SS^{q-1}|\sum_{k=1}^N|H_n(\x,\xi_k)|&=&\sum_{k=1}^N\left|\sum_{i=1}^{N}(\lambda_i-\lambda_{i+1})
\sum_{j=1}^iH_n(\x,\xi_j)K_n(\xi_k\cdot\xi_j)\right|
\\
&\le&\sup_{1\le r\le N}\left(\sum_{k=1}^N\left|\sum_{i=r}^{N}(\lambda_i-\lambda_{i+1})
\sum_{j=1}^iH_n(\x,\xi_j)K_n(\xi_k\cdot\xi_j)\right|\right)
\\
&\le&\sup_{1\le r\le N}\left(\sum_{k=1}^N\sum_{i=r}^{N}(\lambda_i-\lambda_{i+1})
\left|\sum_{j=1}^iH_n(\x,\xi_j)K_n(\xi_k\cdot\xi_j)\right|\right)
\\
&=&\sup_{1\le r\le N}\left(\sum_{i=r}^{N}(\lambda_i-\lambda_{i+1})\sum_{k=1}^N\left|
\sum_{j=1}^iH_n(\x,\xi_j)K_n(\xi_k\cdot\xi_j)\right|\right)
\\
&\le&\sup_{1\le r\le N}\left(\sum_{i=r}^{N}(\lambda_i-\lambda_{i+1})\sup_{r\le i\le N}\sum_{k=1}^N\left|
\sum_{j=1}^iH_n(\x,\xi_j)K_n(\xi_k\cdot\xi_j)\right|\right)
\\
&=&\sup_{1\le r\le N}\left(\lambda_r\sup_{r\le i\le N}\sum_{k=1}^N\left|
\sum_{j=1}^iH_n(\x,\xi_j)K_n(\xi_k\cdot\xi_j)\right|\right).
\end{eqnarray*}
Hence, set
\[
A_r(\x):=\lambda_r\sup_{r\le i\le N}\sum_{k=1}^N\left|
\sum_{j=1}^iH_n(\x,\xi_j)K_n(\xi_k\cdot\xi_j)\right|, \qquad r=1,\ldots,N,
\]
we are going to prove that as $n\rightarrow +\infty$, we have
\begin{equation}\label{eq-th-tens}
\sup_{x\in X_N}A_r(\x)=O( n^{\frac{q-1}2}),\qquad  r=1,\ldots,N .
\end{equation}

First of all, we prove (\ref{eq-th-tens}) for $r=N$. Indeed, from (\ref{inva}) and $\lambda_N=\min_{1\le k\le N} \lambda_k$, we deduce that for all $\x\in X_N$
\[
A_N(\x):=\sum_{k=1}^N\lambda_{N}\left|\sum_{j=1}^NH_n(\x,\xi_j)
K_n(\xi_k\cdot\xi_j)\right|=\sum_{k=1}^N\lambda_{N}|K_n(\x\cdot\xi_k)|\le \sum_{k=1}^N\lambda_{k}|K_n(\x\cdot\xi_k)|.
\]
On the other hand,  it is known \cite{W-Dai, W-Re} that the following Marcinkiewicz inequality follows from the existence of the quadrature rule (\ref{quad-2n})
\begin{equation}\label{Marci}
\sum_{i=1}^N \lambda_i |Q(\xi_i)|\le \C \|Q\|_{L^1(\SS^q)}, \qquad \forall Q\in\PP_{n}, \qquad \C\ne\C(n,N,Q).
\end{equation}
Hence, by using (\ref{Marci}), (\ref{Leb-Fou}) and (\ref{min-norm}), the previous estimate continues as follows
\[
A_N(\x)\le \sum_{k=1}^N\lambda_{k}|K_n(\x\cdot\xi_k)|
\le\C\int_{\SS^q}|K_n(\x\cdot\y)|d\sigma(\y)\le\C\|\S_n\|_\infty \le \C n^{\frac{q-1}2},
\]
i.e.,  we get
\begin{equation}\label{AN-1}
\sup_{\x\in X_N} A_N(\x)\le \C n^{\frac{q-1}2}, \qquad \C\ne\C(n,N).
\end{equation}
Now, for any $1\le r< N$ let us prove that
the assumption
\begin{equation}\label{hp}
\lambda_{r+1}\le \lambda_r\le\C\lambda_{r+1}, \qquad \C\ne\C(n,N,r),
\end{equation}
implies
\begin{equation}\label{th-ind}
\sup_{\x\in X_N} A_{r+1}(\x)\le  \sup_{\x\in X_N} A_r(\x)\le  2 \C \sup_{\x\in X_N} A_{r+1}(\x),
\end{equation}
where the constant $\C$ in (\ref{th-ind}) is the same of that in (\ref{hp}).

Indeed for any $\x\in X_N$, by the first inequality in (\ref{hp}), we get
\begin{eqnarray*}
A_{r+1}(\x)&:=&  \lambda_{r+1}\sup_{r+1\le i\le N}\sum_{k=1}^N\left|
\sum_{j=1}^iH_n(\x,\xi_j)K_n(\xi_k\cdot\xi_j)\right|\\
 &\le& \lambda_{r}\sup_{r+1\le i\le N}\sum_{k=1}^N\left|
\sum_{j=1}^iH_n(\x,\xi_j)K_n(\xi_k\cdot\xi_j)\right|\\
&\le& \lambda_{r}\sup_{r\le i\le N}\sum_{k=1}^N\left|\sum_{j=1}^iH_n(\x,\xi_j)K_n(\xi_k\cdot\xi_j)\right|=  A_{r}(\x),
\end{eqnarray*}
which yields the first inequality in (\ref{th-ind}).

In order to state the second inequality in (\ref{th-ind}), we distinguish two cases.

{\it Case 1:} $\displaystyle \sup_{r+1\le i\le N}\sum_{k=1}^N\left|
\sum_{j=1}^iH_n(\x,\xi_j)K_n(\xi_k\cdot\xi_j)\right|= \sup_{r\le i\le N}\sum_{k=1}^N\left|
\sum_{j=1}^iH_n(\x,\xi_j)K_n(\xi_k\cdot\xi_j)\right|$.

In this case, by the second inequality in (\ref{hp}), we get
\begin{eqnarray*}
A_{r}(\x)&:=& \lambda_{r} \sup_{r\le i\le N}\sum_{k=1}^N\left|
\sum_{j=1}^iH_n(\x,\xi_j)K_n(\xi_k\cdot\xi_j)\right|\\
&=& \lambda_{r} \sup_{r+1\le i\le N}\sum_{k=1}^N\left|
\sum_{j=1}^iH_n(\x,\xi_j)K_n(\xi_k\cdot\xi_j)\right|\\
&\le& \C \lambda_{r+1} \sup_{r+1\le i\le N}\sum_{k=1}^N\left|
\sum_{j=1}^iH_n(\x,\xi_j)K_n(\xi_k\cdot\xi_j)\right|=\C A_{r+1}(\x).
\end{eqnarray*}
{\it Case 2:} $\displaystyle \sup_{r+1\le i\le N}\sum_{k=1}^N\left|
\sum_{j=1}^iH_n(\x,\xi_j)K_n(\xi_k\cdot\xi_j)\right|< \sup_{r\le i\le N}\sum_{k=1}^N\left|
\sum_{j=1}^iH_n(\x,\xi_j)K_n(\xi_k\cdot\xi_j)\right|$.

In this case, by taking into account that
\begin{eqnarray*}
 &&\sup_{r\le i\le N}\sum_{k=1}^N\left|
\sum_{j=1}^iH_n(\x,\xi_j)K_n(\xi_k\cdot\xi_j)\right|\\
&=& \max\left\{\sup_{r+1\le i\le N}\sum_{k=1}^N\left|
\sum_{j=1}^iH_n(\x,\xi_j)K_n(\xi_k\cdot\xi_j)\right|, \ \sum_{k=1}^N\left|
\sum_{j=1}^{r}H_n(\x,\xi_j)K_n(\xi_k\cdot\xi_j)\right|\right\},
\end{eqnarray*}
we can say that
\[
 \sup_{r\le i\le N}\sum_{k=1}^N\left|
\sum_{j=1}^iH_n(\x,\xi_j)K_n(\xi_k\cdot\xi_j)\right|= \sum_{k=1}^N\left|
\sum_{j=1}^{r}H_n(\x,\xi_j)K_n(\xi_k\cdot\xi_j)\right|.
\]
Consequently,  by the second inequality in (\ref{hp}), we get
\begin{eqnarray*}
A_{r}(\x)&=&\lambda_{r}\sum_{k=1}^N\left|
\sum_{j=1}^rH_n(\x,\xi_j)K_n(\xi_k\cdot\xi_j)\right|\\
&\le & \C \lambda_{r+1}\sum_{k=1}^N\left|
\sum_{j=1}^{r+1}H_n(\x,\xi_j)K_n(\xi_k\cdot\xi_j)- H_n(\x,\xi_{r+1})K_n(\xi_k\cdot\xi_{r+1})\right|\\
&\le& \C \lambda_{r+1}\sum_{k=1}^N\left|
\sum_{j=1}^{r+1}H_n(\x,\xi_j)K_n(\xi_k\cdot\xi_j)\right|+
\C\lambda_{r+1}\sum_{k=1}^N\left|
H_n(\x,\xi_{r+1})K_n(\xi_k\cdot\xi_{r+1})\right|\\
&=& \C A_{r+1}(\x) + \C \lambda_{r+1} \sum_{k=1}^N\left|
H_n(\x,\xi_{r+1})K_n(\xi_k\cdot\xi_{r+1})\right|.
\end{eqnarray*}
Moreover, by means of (\ref{ls-ker1}) and (\ref{inva}), we observe that
\begin{eqnarray*}
\lambda_{r+1} \sum_{k=1}^N\left|H_n(\x,\xi_{r+1})K_n(\xi_k\cdot\xi_{r+1})\right|&\le&
\lambda_{r+1} \sum_{k=1}^N\left|K_n(\xi_k\cdot\xi_{r+1})\right|\\
&=& \lambda_{r+1} \sum_{k=1}^N\left|\sum_{j=1}^N H_n(\xi_{r+1}, \xi_j)K_n(\xi_k\cdot\xi_j)\right|\\
&\le& \lambda_{r+1}\sup_{r+1\le i\le N}\sum_{k=1}^N\left|
\sum_{j=1}^iH_n(\xi_{r+1}, \xi_j)K_n(\xi_k\cdot\xi_j)\right|\\
&=& A_{r+1}(\xi_{r+1}).
\end{eqnarray*}
Hence, in the second case we conclude that
\[
A_r(\x)\le \C A_{r+1}(\x)+\C A_{r+1}(\xi_{r+1}).
\]
Summing up, in both the previous cases, for all $\x\in X_N$, we can say that
\[
A_{r}(\x)\le \C A_{r+1}(\x)+ \C A_{r+1}(\xi_{r+1})\le 2\C\sup_{\xi\in X_N} A_{r+1}(\xi),
\]
which yields the second inequality in (\ref{th-ind}).

In conclusion, let us prove that (\ref{th-ind}) and (\ref{AN-1}) imply (\ref{eq-th-tens}).

Indeed,
set for brevity
\[
A_r:=\sup_{\x\in X_N} A_r(\x),\qquad r=1,\ldots,N,
\]
we have to prove that
\[
\limsup_{n\rightarrow +\infty}\ n^{\frac{1-q}2}A_r < +\infty,\qquad r=1,\ldots,N.
\]
But if ad absurdum, for some index $l$ we have
\begin{equation}\label{hp-abs}
\limsup_{n\rightarrow +\infty}\ n^{\frac{1-q}2}A_l= +\infty,
\end{equation}
then by virtue of (\ref{AN-1}) it will be $1\le l <N$. Consequently, since (\ref{th-ind}) yields
\[
\frac{1}{2\C}\ n^{\frac{1-q}2} A_l\le n^{\frac{1-q}2} A_{l+1}\le n^{\frac{1-q}2} A_l,\qquad \C\ne \C(n,N,l),
\]
from (\ref{hp-abs}) we deduce that
\[
\limsup_{n\rightarrow +\infty}\ n^{\frac{1-q}2}A_{l+1}= +\infty
\]
holds too. Then by iterating the reasoning, we arrive to say that
\[
\limsup_{n\rightarrow +\infty}\ n^{\frac{1-q}2} A_N = +\infty,
\]
which contradicts (\ref{AN-1}).

Hence, due to (\ref{th-ind}) and (\ref{AN-1}), we conclude that (\ref{eq-th-tens}) necessarily holds.
\section{Conclusions}
On the unit sphere $\SS^q\subset\RR^{q+1}$, with $q\ge 2$ we studied the approximation provided by least squares polynomials, $\tilde\S_nf$ defined in (\ref{LS-min}), w.r.t. the uniform norm. To this aim, we estimated the behaviour of the associated Lebesgue constants as the polynomial degree $n$ tends to infinity.

Similarly to the hyperinterpolation approximation, for all the polynomial degrees $n$, we supposed that the underlying point set $X_N=\{\xi_1,\ldots,\xi_N\}$ consists of nodes of a positive weighted quadrature rule of degree of precision $2n$.

Then, for least squares polynomial approximation, we stated an optimal behaviour of Lebesgue constants by proving that they grow at the minimal projection order (namely as $n^{\frac{q-1}2}$) under two different additional hypotheses:
\begin{itemize}
\item In a first case (cf.\ Theorems \ref{th-LS} and \ref{th-LSequi}) we supposed that the Marcinkiewicz type inequality (\ref{Marci-q}) holds. This is equivalent to requiring that the nodes $\xi_i$ in $X_N$ are well--separated, namely (\ref{hp-sep}) holds.
\item In the second case (cf. Theorem \ref{th-tensor}), the nodes  can be also not well--separated, but we require that the quadrature weights $\lambda_i$, labeled in non increasing order, satisfy (\ref{hp-tensor}).
\end{itemize}
We remark that in the literature one can find a variety of quadrature nodes fitting into the first or the second case (see, e.g., \cite{W4, r939, W-Po, W-SloBook, W-Xu1} ). In particular, a point set satisfying Theorem \ref{th-LSequi} can be selected from any sufficiently dense set of points on the sphere \cite{W2, W-GiaMha, W-NPW}. Moreover, Theorem \ref{th-tensor} can be applied to the tensor product Gauss--Legendre nodes in \cite{W4}.

In conclusion, under our assumptions, we can say that the approximation provided by least squares and hyperinterpolation polynomials are comparable w.r.t.\ the uniform norm, having in both cases optimal Lebesgue constants.

From a computational point of view, least squares polynomials depend only on the function values at the nodes, while hyperinterpolation polynomials also require a preliminary knowledge of the quadrature weights. Hence, the choice of hyperinterpolation or least squares polynomial approximation depends on the specific problem at hand.

\bibliographystyle{abbrv}
\bibliography{longstrings,TOTAL}
\end{document}